\documentstyle[amscd, amstex]{amsart}

\def\ZZ         {{\Bbb Z}}
\def\RR         {{\Bbb R}}
\def\CC         {{\Bbb C}}

\def\PP         {{\Bbb P}}

\def\AA         {{\Bbb A}}

\newtheorem{thm}{Theorem}[section]
\newtheorem{lem}[thm]{Lemma}
\newtheorem{cor}[thm]{Corollary}
\newtheorem{pr}[thm]{Proposition}
\newtheorem{conj}[thm]{Conjecture}
\theoremstyle{definition}
\newtheorem{rem}[thm]{Remark}
\newtheorem{ex}[thm]{Example}

\newcommand{\key}{\bibitem}

\newcommand{\ps}{{{\bf P}_{\Sigma}}}

\newcommand\hidot{{\raise1pt\hbox{$\scriptscriptstyle\bullet$}}}
\newcommand\lodot{{\raise.3pt\hbox{$\scriptscriptstyle\bullet$}}}

\newcommand{\pp}{{\bf P}}

\begin{document}

\title[Embedding of Calabi-Yau deformations]
{Embedding of Calabi-Yau deformations \\ into toric varieties}
\author{Anvar R. Mavlyutov}
\address {Department of Mathematics, Indiana University, Bloomington, IN 47405, USA.}
 \email{amavlyut@@indiana.edu}


\maketitle

\tableofcontents

\section{Introduction.}

In \cite{m3}, we have constructed non-polynomial deformations of
Calabi-Yau hypersurfaces in toric varieties. Here, we will
construct an embedding of such deformations into toric varieties.
More precisely, we show that deformations of semiample (minimal)
Calabi-Yau hypersurfaces in complete simplicial toric varieties
can be realized as Calabi-Yau complete intersections in higher
dimensional simplicial toric varieties. This result can be
extended to all nef Calabi-Yau complete intersections as well.

There are several questions and conjectures which arise as a
result of this construction. The deformations of semiample
Calabi-Yau hypersurfaces (and, more generally, nef complete
intersections) are realized as non-nef quasi-smooth complete
intersections in higher dimensional toric varieties (with the sum
of the degrees of the hypersurfaces  anticanonical).
 Such complete intersections have not been studied so far.
It will be nice to have a similar deformation result for them.
Because, in this case, it will show that Calabi-Yau complete
intersections in toric varieties preserve the property of being
complete intersections in toric varieties under complex
deformations. As an application to the classification of
Calabi-Yau varieties, it will suffice to look for those Calabi-Yau
complete intersections in toric varieties that have only
``polynomial'' deformations, which means that all deformations are
realized inside the ambient toric variety. Certainly, we expect
the result to hold for singular Calabi-Yau complete intersections
as well.
 This also leads us to the following
question: is it true that every Calabi-Yau variety can be realized
as a   complete intersection in a toric variety?
While this may not look believable at the moment with vastly different constructions of Calabi-Yau varieties (e.g., complete
intersections
in Grassmanians and Flag Manifolds, Voisin-Borcea construction and the Pfaffian Calabi-Yau manifold),
we want to make a prediction of a seemingly achievable goal:
\begin{conj} Extremal transitions (contractions followed by smoothing)
of Calabi-Yau complete intersections in toric varieties are again Calabi-Yau complete
intersections in toric varieties and such Calabi-Yau varieties are all connected to each other by chains of extremal transitions.
\end{conj}

In this conjecture, we assume that the complete intersections of
hypersurfaces have the property that the sum of the degrees is
anticanonical. In physics, it was expected that the extremal
transitions  of Calabi-Yau varieties connect them all together
into a single web. As a consequence of the above conjecture,
Calabi-Yau complete intersections in Grassmanians and Flag
Manifolds must be realized as complete intersections in toric
varieties, because they degenerate to Calabi-Yau complete
intersections in toric varieties. We  also want to make another
conjecture that  Grassmanians and Flag manifolds can be realized
as complete intersections in toric varieties. The reason for this
is that the deformations of Calabi-Yau hypersurfaces in toric
varieties are induced by the deformations of the ambient toric
varieties, which have been realized as complete intersections in
another toric variety. But it is known that Grassmanians and Flag
manifolds are deformations of toric varieties. In particular,
describing deformations of (singular) toric varieties is  an
important problem.

The first example of non-polynomial deformations of a Calabi-Yau
hypersurface in  $\PP(1,1,2,2,2)$, due to Sheldon Katz and David
Morrison (see \cite{cdfkm}), was realized as a $(4,2)$ complete
intersection in $\PP^5$. Even though this example was helpful in
creating our general construction, it turned out that there is
some flexibility in realizing the deformations of Calabi-Yau
hypersurfaces as complete intersections in toric varieties.
However, it didn't help to avoid the following technical problem.
We cannot make all deformations of a Calabi-Yau hypersurface $X$
simultaneously as complete intersections in one toric variety, but
rather we do combine the polynomial deformations together with all
those non-polynomial deformations which correspond, by the
Kodaira-Spencer map,  to a basis in $H^1(X,{\cal T}_X)$ associated
to one of the edges of the reflexive polytope.  Despite this
technical difficulty, the next  priority seems to be in studying
the new Calabi-Yau complete intersections achieved in this paper.
Since, these complete intersections are non-nef, the cohomology of
them will be harder to compute than the cohomology of  the nef
ones, which in the case of nef nondegenerate intersections follows
from \cite{m1+} using a ``Cayley'' trick. Fortunately, the
complete intersections in our construction are quasismooth (i.e.,
V-submanifolds) in simplicial toric varieties, which should help
to compute its cohomology.

The plan of our paper is to quickly  introduce the reader to the
toric varieties and review the results on the infinitesimal
deformations of Calabi-Yau hypersurfaces in Section~\ref{s:t}.
Then, in Section~\ref{s:c}, we consider different constructions of
``non-polynomial'' complex deformations of Calabi-Yau
hypersurfaces which are crepant resolutions of Calabi-Yau
hypersurfaces in weighted projective spaces. The last section
contains our general construction which realizes the ambient toric
variety, the Calabi-Yau hypersurface in it and the non-polynomial
deformations  as quasismooth complete intersections in a higher
dimensional simplicial toric variety.

\section{Toric varieties, Calabi-Yau hypersurfaces and
deformations.}\label{s:t}

This section  reviews  the definition of a  toric variety, and
some facts about  Calabi-Yau hypersurfaces. References on toric
varieties include \cite{c,f,d,o}.

As usual, let $N$ be a lattice
 of rank $d$. A  fan $\Sigma$ is
 a finite collection of strongly convex rational
polyhedral cones in $N_\RR$ such that
   each face of a cone in $\Sigma$
belongs to $\Sigma$ and the intersection of two cones in $\Sigma$
is a face of each.

Assume that 1-dimensional cones of $\Sigma$ span $N_\RR$. Let
$x_1,\dots,x_n$ be coordinates on $\CC^n$, corresponding to the
minimal integral generators $e_1,\dots,e_n$ of the 1-dimensional
cones in $\Sigma$, and let $B(\Sigma)=\langle
\prod_{e_i\not\in\sigma} x_i: \sigma\in\Sigma\rangle$ be the ideal
in $\CC[x_1,\dots,x_n]$. Consider the complement to the closed
subset defined by $B(\Sigma)$, which  is a Zariski open set
$\CC^n\setminus{\bf V}(B(\Sigma))$. This set is  invariant under
the diagonal action of an affine algebraic group $G=G(\Sigma)$, a
subgroup of $(\CC^*)^n$, defined by the  equations $\prod_{i=1}^n
x_i^{\langle m,e_i\rangle}=1$, for  $m\in M={\rm Hom}(N,\ZZ)$,
where $\langle\,,\,\rangle$ is the pairing. The toric variety
$\ps$, associated with the fan $\Sigma$, is the categorical
quotient  $(\CC^n\setminus{\bf V}(B(\Sigma))/G$. The ring
$S:=S(\Sigma)={\Bbb C}[x_1,\dots,x_n]$ is   the homogeneous
coordinate ring
 of $\ps$ graded by the Chow group
$A_{d-1}(\ps)$:  $\deg(\prod_{i=1}^n x_i^{a_i})=[\sum_{i=1}^n a_i
D_i]$, where $D_i$ is an irreducible divisor given by $x_i=0$.

 Semiample anticanonical nondegenerate
(transversal to the torus orbits) hypersurfaces  in a complete
simplicial toric varieties  are important examples of Calabi-Yau
varieties.  They are defined by the homogeneous polynomials $f\in
S_\beta$ (here, $\beta=[\sum_{i=1}^n  D_i]$ is the anticanonical
degree) such that $x_i\partial f/\partial x_i$ do not vanish
simultaneously on $\ps$. The fact that $D=\sum_{i=1}^n D_i$ is
semiample, implies that the associated polytope
$$\Delta:=\Delta_D=\{m\in M_{\Bbb R}:\langle m,e_i\rangle\geq-1
\text{ for all } i\}\subset M_{\Bbb R}$$ is reflexive, which means
all of its vertices are lattice points, the only interior point is
the origin $0$, and the same property has to hold for the dual
polytope $$\Delta^\circ=\{n\in N_{\Bbb R}:\langle m,n\rangle\geq-1
\text{ for all } m\in\Delta\}\subset M_{\Bbb R}.$$ Moreover, the
minimal integral generators $e_1,\dots,e_n$ of the 1-dimensional
cones lie on the boundary of $\Delta^\circ$. Since
 $X$ is a big and nef hypersurface, by Proposition 1.2 in \cite{m1},
 there is the associated toric
morphism $\pi:{\bf P}_\Sigma\rightarrow{\bf P}_{\Sigma_X}$, where
$\Sigma_X$ is the fan whose cones are precisely those  generated
by the faces of the polytope $\Delta^\circ$. For instance, a
2-dimensional cone in $\Sigma_X$ is a cone over an edge
$\Gamma^*\subset \Delta^\circ$.

 When a variety $X$ is smooth, its space of
infinitesimal deformations is classified by the first cohomology
$H^1(X,{\cal T}_X)$ with coefficients in the tangent sheaf ${\cal
T}_X$. For semiample anticanonical nondegenerate hypersurfaces
$X\subset\ps$ in a complete simplicial toric variety, there is a
natural splitting $$H^1(X,{\cal T}_X)=H_{\rm poly}^1(X,{\cal
T}_X)\oplus H_{\rm nonp}^1(X,{\cal T}_X),$$ where $H_{\rm
poly}^1(X,{\cal T}_X)$ is the image of the map $H^0(X,{\cal
N}_{X/\ps})@>>>H^1(X,{\cal T}_X)$, which classifies all
deformations of $X$ inside $\ps$. While the space $H_{\rm
nonp}^1(X,{\cal T}_X)$ of non-polynomial infinitesimal
deformations of $X$ can be thought as the cokernel of the last
map, a natural description is given in Theorem~7.1 in \cite{m2}
(also, see \cite[Section 1]{m3}).

\begin{pr}\label{t:chiral}
Let $X\subset\ps$ be  a  semiample anticanonical nondegenerate
 hypersurface in a complete
simplicial toric variety, defined by $f\in S_\beta$. Then there
are natural  isomorphisms $$ H^1_{\rm poly}(X,{\cal T}_X)\cong
R_1(f)_{\beta},\qquad H^1_{\rm nonp}(X,{\cal
T}_X)\cong\bigoplus\begin{Sb}\sigma \in\Sigma_X,\,\dim\sigma=2\\
e_l\in{{\rm int}(\sigma)}\end{Sb}(S/\langle
x_l\rangle)_{\beta_1^\sigma},$$ where $R_1(f)=S/J_1(f)$ with
$J_1(f)=\langle x_1(\partial f/\partial x_1),\dots,x_n(\partial
f/\partial x_n)\rangle:x_1\cdots x_n$, and
$\beta^\sigma_1=\sum_{e_i\subset\sigma}\deg(x_i)$.
\end{pr}

\begin{cor}\label{p:dim} Let $X$ be a semiample anticanonical nondegenerate hypersurface
in a complete simplicial toric variety $\ps$, such that the
minimal lattice generators $e_1,\dots,e_n$ of the 1-dimensional
cones in $\Sigma$ include all lattice points in the edges of
$\Delta^\circ$, the dual of the reflexive polytope $\Delta$
associated to the anticanonical divisor.
 Then
 $$\dim
H_{\rm poly}^1(X,{\cal T}_X)=l(\Delta)-d-1-\sum_{\rm codim
\Gamma=1}l^*(\Gamma),$$ $$\dim H_{\rm nonp}^1(X,{\cal
T}_X)=\sum_{\rm codim \Gamma=2}l^*(\Gamma)l^*(\Gamma^*), $$
 where the sum is by faces $\Gamma$ of $\Delta$, $\Gamma^*$ is the corresponding
  dual face of $\Delta^*$
  and
$l(\Delta)={\rm Card}(\Delta\cap M)$, while $l^*(\Gamma)$ is the
number of interior lattice points in $\Gamma$. Moreover, if
$e_1,\dots,e_n$ contain $\Delta^\circ\cap N-\{0\}$, except,
possibly, the interior points of  the facets of $\Delta^\circ$,
then $X$ is a minimal Calabi-Yau orbifold.
\end{cor}

\begin{pf} The fact that the hypersurface $X$ is a minimal Calabi-Yau
orbifold is essentially shown in Proposition 4.1.3 of \cite{ck}.
One simply has to note that the singularities of a semiample
nondegenerate hypersurface come from singularities of the toric
variety, which are all terminal by \cite[Lemma 4.1.2]{ck}. Hence,
singularities of $X$ are also terminal.

 Arguments given in
Remark 5.3 of \cite{m1} combined with the description of the space
$H^1(X,{\cal T}_X)$ in Proposition \ref{t:chiral} produce the
required dimension formulas.
\end{pf}

We can also determine the polynomial defining a partial crepant
resolution of a given Calabi-Yau hypersurface $Y$ in a complete
toric variety $\pp_{\Sigma'}$. Let $Y$ be determined by a
homogeneous polynomial
 $$\sum_{m\in\Delta\cap M}a_m \prod_{i=1}^n
z_i^{1+\langle m,e'_i\rangle}$$ in $S(\Sigma')$, where $\Delta$ is
the polytope in $M_\RR$ given by the conditions $\langle
m,e'_i\rangle\ge-1$ and $e_i'$ are lattice generators of the
1-dimensional cones in $\Sigma'$. Then, if the fan $\Sigma$ is a
subdivision of $\Sigma'$, the proper transform $X\subset\ps$ of
the hypersurface $Y$ corresponds to $$\sum_{m\in\Delta\cap M}a_m
\prod_{i=1}^n x_i^{1+\langle m,e_i\rangle}$$ in $S(\Sigma)$.

\section{Calabi-Yau deformations arising from weighted projective
spaces.}\label{s:c}

First, we will  show how our construction works for some examples.

The simplest example of a non-polynomial deformation of a
Calabi-Yau hypersurface  arises from the weighted projective space
$\PP(1,1,2,2,2)$. This case, considered in \cite{cdfkm}, was
described in details in \cite[Example 2.1]{m3}. Let us recall this
construction.

\begin{ex}\label{e:p11}
 Consider the quasismooth Calabi-Yau hypersurface $Y$
in $\PP(1,1,2,2,2)$
 defined  by the degree 8 polynomial
 $$z_1^ 8 + z_2 ^8 + z_3^ 4
+ z_4^ 4 + z_5^ 4 = 0.$$

The crepant desingularization of $\PP(1,1,2,2,2)$ is the toric
variety
 $\ps$ whose fan has the  following integral generators of
 the 1-dimensional cones:
\begin{equation}\label{e:fan}
\{(-1,-2,-2,-2),(1,0,0,0),(0,1,0,0),(0,0,1,0),(0,0,0,1),(0,-1,-1,-1)\}.
\end{equation}
If $x_1,x_2,x_3,x_4,x_5,x_6$ denote the homogeneous coordinates of
$\ps$, corresponding to these points, then the crepant resolution
$X$ of $Y$ is defined by the polynomial $$f(x)=x_1^ 8x_6^4 + x_2
^8x_6^4 + x_3^4+ x_4^ 4 + x_5^ 4.$$

Both the toric variety $\ps$ and the minimal Calabi-Yau
hypersurface $X$ can be embedded into the product
$\PP^5\times\PP^1$ by sending
\begin{equation*}
(x_1,x_2,x_3,x_4,x_5,x_6)\mapsto(x_1x_2x_6,x_1^2x_6,x_2^2x_6,x_3,x_4,x_5)\times(x_1,x_2).
\end{equation*}

It is straightforward to check that this is an inclusion. If
$y_0,\dots,y_5$ and $z_1,z_2$ are the homogeneous coordinates on
the 2 factors, the hypersurface $X$ becomes a Calabi-Yau complete
intersection given by $$y_1^4 + y_2 ^4 + y_3^ 4 + y_4^ 4 + y_5^ 4
= 0, \qquad y_2z_1=y_0z_2,\qquad y_1z_2=y_0z_1,$$ where the last
two equations determine the variety $\ps$.

Using the dimension formulas in Proposition~\ref{p:dim} and a
computer program (e.g., see \cite{ks}), one can find that $\dim
H_{\rm poly}^1(X,{\cal T}_X)=83$ and $\dim H_{\rm nonp}^1(X,{\cal
T}_X)=3$.
 Polynomial deformations of the hypersurface
$X\subset\ps$, which leave the hypersurface inside its ambient
space $\ps$, are created by deforming the coefficients of the
quartic, while the rest --- the non-polynomial deformations
--- correspond to the family $$y_2z_1=(y_0+\lambda_3 y_3+\lambda_4
y_4 +\lambda_5 y_5 )z_2,\quad y_1z_2=y_0z_1,$$ which also induce
 deformations of the toric variety $\ps$.

 The projection of the non-polynomial deformations of $X$  onto the first factor
$\PP^5\times\PP^1\rightarrow \PP^5$ gives the Calabi-Yau complete
intersection
  $$y_1^4 + y_2 ^4 + y_3^ 4 + y_4^ 4 + y_5^ 4 = 0,\qquad y_1y_2=y_0(y_0+\lambda_3 y_3+\lambda_4
y_4 +\lambda_5 y_5),$$
 which is a deformation of the quasismooth hypersurface $Y$ in
 $\PP(1,1,2,2,2)$: the hypersurface  $Y$ together with its ambient space
 are embedded into $\PP^5$ by  the degree 2 linear system:
 $$(z_1,z_2,z_3,z_4,z_5)\mapsto(z_1z_2,z_1^2,z_2^2,z_3,z_4,z_5).$$
   Moreover, the  Calabi-Yau variety $Y$
  can be smoothed to a generic complete intersection of
degree 4 and 2 hypersurfaces, and the process of going from $X$ to
this smooth Calabi-Yau complete intersection is called an extremal
transition (see \cite[Definition 6.2.9.]{ck}).

We now make an important observation that the deformation of $Y$
is realized as a complete intersection of 2 hypersurfaces, while
$X$ requires 3. This leads us to the question whether the
deformations of $X$ can also be realized as a complete
intersection of 2 hypersurfaces, which we are going to answer
positively.

Embed the hypersurface $X$ with its ambient space $\ps$ into
 $\PP^4\times\PP^1$ using the map
\begin{equation}\label{e:map1}
(x_1,x_2,x_3,x_4,x_5,x_6)\mapsto(x_1^2x_6,x_2^2x_6,x_3,x_4,x_5)\times(x_1,x_2).
\end{equation}
The image of $X$ is a Calabi-Yau complete intersection:
\begin{equation*}\label{e:ci}
y_1^4 + y_2 ^4 + y_3^ 4 + y_4^ 4 + y_5^ 4 = 0,\qquad
y_1z_2^2=y_2z_1^2,
\end{equation*}
 and the last equation cuts out the subvariety
isomorphic to $\ps$ in $\PP^4\times\PP^1$. Then all deformations
of $X$ are given by deforming the degree $(4,0)$ hypersurface and
by the family
\begin{equation*}\label{e:nd} y_1z_2^2=y_2z_1^2+(\lambda_3 y_3+\lambda_4 y_4
+\lambda_5 y_5 )z_1z_2
\end{equation*}
of degree $(1,2)$ hypersurfaces, which correspond to deformations
of the toric variety $\ps$.
\end{ex}

The second example,  we consider in this section, is coming from
$\PP(1,1,2,2,6)$.

\begin{ex} The weighted projective space $\PP(1,1,2,2,6)$
  has  the fan consisting of cones
with generators
\begin{equation}\label{e:fan26}
\{(-1,-2,-2,-6),(1,0,0,0),(0,1,0,0),(0,0,1,0),(0,0,0,1),\}.
\end{equation}
The convex hull of these points is a reflexive polytope, such that
the only extra lattice points on its boundary are $(0,-1,-1,-3)$
and $(0,0,0,-1)$. Note that subdividing the fan using
$(0,-1,-1,-3)$ is enough to get a complete simplicial toric
variety $\ps$ such that the proper transform of a nondegenerate
Calabi-Yau hypersurface $Y$ in $\PP(1,1,2,2,6)$ is smooth. This is
so, since $(0,0,0,-1)$ is in the interior of a 4-dimensional cone,
and the nondegenerate hypersurface $Y$ will miss the corresponding
singular point of  $\PP(1,1,2,2,6)$. Let $Y$ be given by the
degree 12 polynomial $$z_1^{12} + z_2 ^{12} + z_3^ 6 + z_4^ 6 +
z_5^ 2 = 0,$$ where $z_1,\dots,z_5$ correspond to the lattice
points (\ref{e:fan26}). Then its crepant  resolution $X\subset\ps$
is an anticanonical nondegenerate hypersurface $$x_1^{12}x_6^6 +
x_2 ^{12}x_6^6 + x_3^ 6 + x_4^ 6 + x_5^ 2 = 0,$$ where
$x_1,\dots,x_5,x_6$ correspond to (\ref{e:fan26}) and
$(0,-1,-1,-3)$.

In this case, we find $\dim H_{\rm poly}^1(X,{\cal T}_X)=126$ and
$\dim H_{\rm nonp}^1(X,{\cal T}_X)=2$.

Embed the toric variety $\ps$ and the smooth Calabi-Yau
hypersurface $X$ into $\PP(1,1,1,1,3)\times\PP^1$ by the same map
as in (\ref{e:map1}). Then $X$ is isomorphic to a Calabi-Yau
complete intersection in $\PP(1,1,1,1,3)\times\PP^1$  given by
$$y_1^6 + y_2 ^6 + y_3^ 6 + y_4^ 6 + y_5^ 2 = 0,\qquad
y_1z_2^2=y_2z_1^2,$$ and the  non-polynomial deformations,
corresponding to $ H_{\rm nonp}^1(X,{\cal T}_X)$, are given by
replacing the second equation with
 $$
 y_1z_2^2=y_2z_1^2+(\lambda_3 y_3+\lambda_4 y_4)z_1z_2.$$
\end{ex}

Using the above two examples, one can make a prediction that if a
minimal semiample Calabi-Yau hypersurface $X$ in a complete
simplicial toric variety $\ps$ is given, then in order to get the
deformations the ambient variety $\ps$ must be embedded into
another toric variety with the map given in homogeneous
coordinates by
$$(x_1,x_2,x_1^2x_6,x_2^2x_6,x_3,x_4,x_5,x_7,\dots),$$ where
$x_1,x_2,x_6$ correspond to lattice points $e_1,e_2,e_6$ such that
$e_6$ is the middle point between $e_1$ and $e_2$. It turns out
that this will work well as long as there is only one lattice
point on the line segment between $e_1$ and $e_2$. In the next
section we will present a general construction, which in  the
particular situation with $\PP(1,1,2,2,2)$ looks as the following.

\begin{ex}
Let  $X$ be   the minimal Calabi-Yau hypersurface defined by
$$x_1^ 8x_6^4 + x_2 ^8x_6^4 + x_3^4+ x_4^ 4 + x_5^ 4=0$$ in the
toric variety $\ps$ as in Example~\ref{e:p11}. Then consider
another  5-dimensional toric variety ${\bf P}_{\Sigma'}$,  whose
 fan $\Sigma'$ has the following minimal integral
generators along the one-dimensional cones:
\begin{multline*}
\{(1,1,1,1,1),(-1,-2,-2,-2,-1),(0,-1,-1,-1,-1),(0,0,0,0,1),\\(0,1,0,0,0),(0,0,1,0,0),
(0,0,0,1,0)\}.\end{multline*}
 The maximal cones of the fan are spanned by one
of the first two generators and four of the last five. By Exercise
in \cite[page~41]{f}, this toric variety is a $\PP^4$ fiber bundle
over $\PP^1$. Moreover, this is a crepant resolution of a Fano
toric variety and the listed generators are the vertices of a
reflexive polytope. The toric variety $\ps$ can be embedded into
${\bf P}_{\Sigma'}$ by the map
$$(x_1,\dots,x_6)\mapsto(x_1,x_2,x_1x_6,x_2x_6,x_3,x_4,x_5),$$
where components are homogeneous coordinates
$s_1,t_1,s_2,t_0,y_3,y_4,y_5$ (notation is made similar to the one
in the next section), corresponding to the above lattice points.
Then the image of $\ps$ is a hypersurface given by
$s_1t_0=s_2t_1$, and the Calabi-Yau hypersurface $X$ is now a
complete intersection of $s_1t_0=s_2t_1$ with
$s_1^4s_2^4+t_0^4t_1^4+y_3^4+y_4^4+y_5^4$. The non-polynomial
deformations are represented by the family of complete
intersections:
$$s_1t_0=s_2t_1+\lambda_3y_3+\lambda_4y_4+\lambda_5y_5\qquad
s_1^4s_2^4+t_0^4t_1^4+y_3^4+y_4^4+y_5^4$$ with three parameters
$\lambda_3$, $\lambda_4$ and $\lambda_5$.
\end{ex}

\section{Embedding of Calabi-Yau  deformations using one
edge.}

First, we consider a simpler construction. Here is the idea of how
we achieved it. Using the examples from the previous section, we
guess what should be the map that embeds the toric variety with
its Calabi-Yau hypersurface into another toric variety.  A problem
is then to find the combinatorial data that would produce a toric
variety compatible with the map.

While we explained how the construction was found, it is more
convenient  to give the answer for the new ambient space first.
Let  $\ps$ be a complete simplicial toric variety with a semiample
anticanonical degree,  such that the minimal lattice generators
$e_1,\dots,e_n$ of the 1-dimensional cones in $\Sigma$ include all
lattice points in one of the  edges
 $\Gamma^*$
of $\Delta^\circ$. We order the latter lattice points
$e_{i_0},e_{i_1},\dots,e_{i_{k+1}}$ as in the following picture,
where the cones of $\Sigma$ subdivide a 2-dimensional cone over
$\Gamma^*$. $$ \setlength{\unitlength}{1cm}
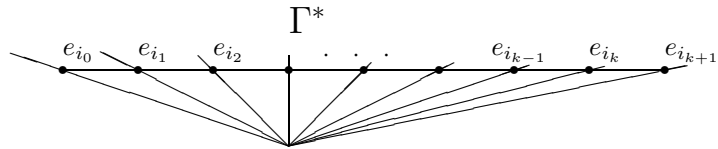
\begin{picture}(8,2)
\put(-1,1.2){$e_{i_0}$} \put(2,0){\line(-3,1){3.7}}
\put(-1,1){\circle*{0.1}} \put(0,1){\circle*{0.1}}
\put(1,1){\circle*{0.1}} \put(2,1){\circle*{0.1}}
\put(3,1){\circle*{0.1}}
\put(4,1){\circle*{0.1}}\put(5,1){\circle*{0.1}}\put(6,1){\circle*{0.1}}\put(7,1){\circle*{0.1}}
 \put(0,1.2){$e_{i_1}$}
\put(2,0){\line(-2,1){2.5}} \put(2,0){\line(-1,1){1.2}}
\put(1,1.2){$e_{i_2}$} \put(2,0){\line(0,1){1.2}}
\put(2,0){\line(1,1){1.1}} \put(2,0){\line(2,1){2.2}}
\put(2,1.5){\LARGE$\Gamma^*$}
\multiput(2.5,1.2)(.4,0){3}{\circle*{0.03}}
\put(4.7,1.2){$e_{i_{k-1}}$} \put(6,1.2){$e_{i_{k}}$}
\put(2,0){\line(3,1){3.2}} \put(2,0){\line(4,1){4.2}}
\put(2,0){\line(5,1){5.3}} \put(7,1.2){$e_{i_{k+1}}$}
\put(-1,1){\line(1,0){8}}
\end{picture}\label{e:pic}
$$

 We want to embed  $\ps$ into another
toric variety, which is  denoted $\pp_{\Sigma(\Gamma^*)}$. The new
toric variety must have a representation as a quotient
$$\pp_{\Sigma(\Gamma^*)}\cong\biggl(\CC^{n+k}\setminus{\bf
V}(B(\Sigma(\Gamma^*))\biggr)/G(\Sigma(\Gamma^*)),$$ for some fan
$\Sigma(\Gamma^*)$. Our goal is to determine this fan.
 Using the
example from the previous section, we first  guess that the
embedding  $\phi_{\Gamma^*}:\ps\rightarrow\pp_{\Sigma(\Gamma^*)}$
can be written   as follows. If $x=(x_1,\dots,x_n)$ represents a
point in $\ps$, then $\phi_{\Gamma^*}(x)$ is the point represented
by   the homogeneous coordinates of $\pp_{\Sigma(\Gamma^*)}$:
\begin{multline}\label{e:map}
s_j=x_{i_0}x_{i_1}\dots x_{i_{j-2}}x_{i_{j-1}},\quad  \text{ for }
j=1,\dots,k+1,  \\  t_j=x_{i_{j+1}}x_{i_{j+2}}\dots
x_{i_k}x_{i_{k+1}},\quad \text{ for } j=0,\dots,k,\text{ and }
y_i=x_i,\text{ for }e_i\notin\Gamma^*,
\end{multline}
where the coordinates $s_j$, $t_j$ correspond to the unknown
lattice points $\tilde e_{s_j}$, $\tilde e_{t_j}$ in the lattice
$N\oplus\ZZ^k$, and $y_i$ correspond to $\tilde e_i=e_i$ in
$N\oplus\ZZ^k$ for $e_i\notin\Gamma^*$. So, the  fan
$\Sigma(\Gamma^*)$ must lie in the lattice $N\oplus \ZZ^{k}$ with
the minimal integral generators of the 1-dimensional cones:
$\tilde e_{s_1},\dots,\tilde e_{s_{k+1}}$, $\tilde
e_{t_0},\dots,\tilde e_{t_k}$ and $\tilde e_i$ for
$e_i\notin\Gamma^*$. Instead of explaining how to find the unknown
lattice points, we will rather give an answer right away. The idea
of deriving these lattice points is to use the fact that the map
$\phi_{\Gamma^*}(x)$ has to be equivariant with respect to the
group actions and it will appear in  the proof of the theorem
below. Let $v_1,\dots,v_k$ be a standard basis of $\ZZ^k$, then
set
\begin{equation}\label{e:generators}
\matrix \tilde e_{s_j}=v_j, \text{ for } j=1,\dots,k,\hfill
 \\ \tilde e_{s_{k+1}}=e_{i_0}-v_1-\cdots-v_k,\hfill
\\
\tilde e_{t_{j-1}}=v_j+\frac{1}{k+1}(e_{i_{k+1}}-e_{i_0}),\text{
for } j=1,\dots,k,\hfill\\ \tilde
e_{t_k}=e_{i_0}-v_1-\cdots-v_k+\frac{1}{k+1}(e_{i_{k+1}}-e_{i_0}).\hfill
\endmatrix
\end{equation}

Next, form  a collection $\Sigma(\Gamma^*)$ of cones in
$N_\RR\oplus \RR^k$, consisting of  the cones with their faces  of
the following types:
\\  ${\rm (i)} \quad \sigma\in\Sigma\quad \text{ if }\quad
 \sigma\cap\Gamma^*=\emptyset, $\\
${\rm (ii)} \quad  \langle e_i\in\sigma\setminus\{
e_{i_j}\},\tilde{e}_{t_0},
\tilde{e}_{t_{1}},\dots,\tilde{e}_{t_{j-1}},\tilde{e}_{s_{j+1}},\dots,\tilde
e_{s{_{k+1}}} \rangle \text{ if }
 \sigma\cap\Gamma^*=\{e_{i_j}\}$,\\
  ${\rm (iii)} \quad  \langle
e_i\in\sigma\setminus(\sigma\cap\Gamma^*),\tilde{e}_{t_0},
\tilde{e}_{t_{1}},\dots,\tilde{e}_{t_{j}},\tilde{e}_{s_{j+1}},\dots,\tilde{e}_{s{_{k+1}}}\rangle
\text{ if }
 \sigma\cap\Gamma^*=\{e_{i_j},{e}_{i_{j+1}}\},$\\
where we indicated the generators of the cones and $e_i$ run by
all possible generators in the complements. To get this collection
of cones note which variables vanish in (\ref{e:map}) if
$x_{i_j}=0$. Notice also that the cones in (i) do not have the
full dimension of the space $N_\RR\oplus \RR^k$. So, the above
cones do not necessarily cover the whole space. Magically, the
collection $\Sigma(\Gamma^*)$ is a fan.

\begin{lem}
The collection of cones $\Sigma(\Gamma^*)$ forms a simplicial fan
in the lattice $N\oplus\ZZ^k$.
\end{lem}

\begin{pf} The fact that $\Sigma(\Gamma^*)$ is simplicial follows
from an easy checking that  the generators of the cones of types
{\rm (ii)-(iii)} are linearly independent and the condition that
$\Sigma$ was simplicial. We need to show that the
  cones in $\Sigma(\Gamma^*)$ meet along their faces, or, in other words, that
  the intersection of the cones lies in the face  generated
  by the common minimal lattice generators.

   Observe the
  property satisfied by the cones  in $\Sigma(\Gamma^*)$: if $\sigma'$ denotes the cone in (ii)-(iii) constructed from
  $\sigma\in\Sigma$, then  the  intersection of $\sigma'$ with the
  subspace $N_\RR$ in $N_\RR\oplus \RR^k$ is exactly $\sigma$.
 Indeed, for $\sigma'$ of type (ii), if $$\alpha_{t_0}\tilde
 e_{t_0}+\cdots+\alpha_{t_{j-1}}\tilde
 e_{t_{j-1}}+\alpha_{s_{j+1}}\tilde
 e_{s_{j+1}}+\cdots+\alpha_{s_{k+1}}\tilde
 e_{s_{k+1}}+\sum_{e_i\in\sigma\setminus\{e_{i_j}\}}\alpha_i e_i\in N_\RR\oplus {\bf 0},$$
where all the coefficients $\alpha\ge0$, then, by comparing the
coefficients at the basis vectors $v_1,\dots,v_k$, we deduce
$$\alpha_{t_0}=\cdots=\alpha_{t_{j-1}}=\alpha_{s_{j+1}}=\cdots\alpha_{s_{k+1}}.$$
But in this case, the above sum is $$\alpha_{s_{k+1}}(\tilde
e_{i_0}+\frac{j}{k+1}(e_{i_{k+1}}-e_{i_0}))+\sum_{e_i\in\sigma\setminus\{e_{i_j}\}}\alpha_i
 e_i=\alpha_{s_{k+1}}e_{i_j}+\sum_{e_i\in\sigma\setminus\{e_{i_j}\}}\alpha_i
 e_i,$$
which belongs to $\sigma$. The other inclusion follows from the
equality $$e_{i_j}=\tilde
 e_{t_0}+\cdots+ \tilde
 e_{t_{j-1}}+ \tilde
 e_{s_{j+1}}+\cdots+ \tilde
 e_{s_{k+1}}.$$
Thus $\sigma'\cap N_\RR\oplus {\bf 0}=\sigma\oplus {\bf 0}$.
Verification of this property for type (iii) cones is similar.

Now, we can see that the intersection of type (i) cone
$\sigma_1\oplus {\bf 0}$ with a cone $\sigma_2'$ of type (ii) or
(iii) is exactly $(\sigma_1\cap\sigma_2)\oplus {\bf 0}$, which is
a face of each of the cones.

Next, let us show that  the intersection of two distinct type (ii)
cones $\sigma_1'$ and $\sigma_2'$ is a face of each. A point lying
in the intersection would produce the following equality:
\begin{multline*}\alpha_{t_0}\tilde
 e_{t_0}+\cdots+\alpha_{t_{j-1}}\tilde
 e_{t_{j-1}}+\alpha_{s_{j+1}}\tilde
 e_{s_{j+1}}+\cdots+\alpha_{s_{k+1}}\tilde
 e_{s_{k+1}}+\sum_{e_i\in\sigma_1\setminus\{e_{i_j}\}}\alpha_i e_i\\=
\beta_{t_0}\tilde
 e_{t_0}+\cdots+\beta_{t_{l-1}}\tilde
 e_{t_{l-1}}+\beta_{s_{l+1}}\tilde
 e_{s_{l+1}}+\cdots+\alpha_{s_{k+1}}\tilde
 e_{s_{k+1}}+\sum_{e_i\in\sigma_2\setminus\{e_{i_l}\}}\alpha_i
 e_i,
 \end{multline*}
where $\sigma_1\cap\Gamma^*=\{e_{i_j}\}$,
$\sigma_2\cap\Gamma^*=\{e_{i_l}\}$ and $j\le l$. Again, comparing
coefficients at $v_1,\dots,v_k$, we get
\begin{align*}
&\alpha_{t_a}-\alpha_{s_{k+1}}=\beta_{t_a}-\beta_{s_{k+1}},\quad
a=0,\dots,j-1,\\
&\alpha_{s_{a+1}}-\alpha_{s_{k+1}}=\beta_{t_a}-\beta_{s_{k+1}},\quad
a=j,\dots,l-1,\\
&\alpha_{s_{a+1}}-\alpha_{s_{k+1}}=\beta_{s_{a+1}}-\beta_{s_{k+1}},\quad
a=l,\dots,k.
\end{align*}
Using this  and simplifying the above equality, we have
\begin{multline*}
\alpha_{s_{k+1}}e_{i_j}+\sum_{e_i\in\sigma_1\setminus\{e_{i_j}\}}\alpha_i
 e_i
 \\=(\beta_{t_j}+\cdots+\beta_{t_{l-1}}-(l-j)\beta_{s_{k+1}})\frac{1}{k+1}(e_{i_{k+1}}-
 e_{i_0})+\beta_{s_{k+1}}e_{i_l}+\sum_{e_i\in\sigma_2\setminus\{e_{i_l}\}}\alpha_i
 e_i.\end{multline*}
If $j=l$, then we are done with this case. If $j<l$, consider
first the case  that both $\alpha_{s_{k+1}}$ and $\beta_{s_{k+1}}$
are not
 zero. Then note that there exists a point $p=(1-
\varepsilon)e_{i_j}+ \varepsilon e_{i_l}$, $0\le\varepsilon\le 1$,
between $e_{i_j}$ and $e_{i_l}$ such that the intersection of the
cones $\langle p,e_i\in\sigma_1\setminus\{e_{i_j}\}\rangle$ and
$\langle p,e_i\in\sigma_2\setminus\{e_{i_l}\}\rangle$ is at most
the cone generated by $p$ and $\sigma_1\cap\sigma_2$. We can
rewrite the last equality as
$$\alpha_{s_{k+1}}e_{i_j}-\varepsilon
(\alpha_{s_{j+1}}+\cdots+\alpha_{s_{l}}-(l-j)\alpha_{s_{k+1}})\frac{1}{k+1}(e_{i_{k+1}}-
 e_{i_0})+
\sum_{e_i\in\sigma_1\setminus\{e_{i_j}\}}\alpha_i
 e_i$$
 $$=(1-\varepsilon)(\beta_{t_j}+\cdots+\beta_{t_{l-1}}-(l-j)\beta_{s_{k+1}})\frac{1}{k+1}(e_{i_{k+1}}-
 e_{i_0})+\beta_{s_{k+1}}e_{i_l}+\sum_{e_i\in\sigma_2\setminus\{e_{i_l}\}}\alpha_i
 e_i,$$
or \begin{multline*} \alpha_{s_{k+1}}p-\varepsilon
(\alpha_{s_{j+1}}+\cdots+\alpha_{s_{l}})\frac{1}{k+1}(e_{i_{k+1}}-
 e_{i_0})+
\sum_{e_i\in\sigma_1\setminus\{e_{i_j}\}}\alpha_i
 e_i
 \\=\beta_{s_{k+1}}p+(1-\varepsilon)(\beta_{t_j}+\cdots+\beta_{t_{l-1}})
 \frac{1}{k+1}(e_{i_{k+1}}-
 e_{i_0})\sum_{e_i\in\sigma_2\setminus\{e_{i_l}\}}\alpha_i
 e_i.\end{multline*}
 Assuming that $\alpha_{s_{k+1}},\beta_{s_{k+1}}\ne0$, the above linear combinations would represent a point
 in the intersection of two convex simplicial cones whose
 generators
 are at a point in the plane passing through the origin and $\Gamma^*$
 and at   $e_i$ in $\sigma_1\setminus\{e_{i_j}\}$ and in
 $\sigma_2\setminus\{e_{i_l}\}$,
respectively. By the choice of $p$, we conclude that
$$\alpha_{s_{j+1}}=\cdots=\alpha_{s_{l}}=\beta_{t_{j}}=\cdots=t_{l-1}.$$
Hence, $\alpha_{s_{k+1}}=\beta_{s_{k+1}}$ as well, and we can
deduce that the intersection $\sigma_1'\cap\sigma_2'$  lies in the
face generated
  by the common minimal lattice generators.
In the case when one of $\alpha_{s_{k+1}}$ and $\beta_{s_{k+1}}$
vanishes we can assume, by symmetry, that $\beta_{s_{k+1}}=0$.
This will lead to an equation where on the left hand side there is
point from $\sigma_1$ and on the right side -- a point from a
facet of $\sigma_2$ shifted away from $\sigma_1$ by the vector
$(\beta_{t_j}+\cdots+\beta_{t_{l-1}})(e_{i_{k+1}}-
 e_{i_0}).$
 This implies $\beta_{t_j}=\cdots=\beta_{t_{l-1}}=\alpha_{s_{k+1}}=0$,
 whence $\alpha_{s_{j+1}}=\cdots=\alpha_{s_{l}}=0$ as well.
 Thus, $\sigma_1'\cap\sigma_2'$ is a face of both  $\sigma_1'$ and
 $\sigma_2'$, when these two cones are of type (ii).
 The other  cases of intersections of two cones of type (ii) and
 (iii) can be shown to hold the required property in a similar manner, and left to the
 reader.
 \end{pf}

Here is our first main result which shows that the toric variety
$\ps$ is realized as a complete intersection by the map
$\phi_{\Gamma^*}$.

\begin{thm}\label{t:emb}
Let $\ps$ be a complete simplicial toric variety with a semiample
anticanonical degree, such that the minimal lattice generators
$e_1,\dots,e_n$ of $\Sigma(1)$ lie on the boundary  of a reflexive
polytope $\Delta^\circ$ and include all lattice points in one of
the  edges
 $\Gamma^*\subset\Delta^\circ$. Then the morphism
 $$\phi_{\Gamma^*}:\ps\rightarrow\pp_{\Sigma(\Gamma^*)}$$
(see (\ref{e:map})) is a well defined  embedding. Moreover, the
image of $\ps$ is a quasismooth complete intersection in
$\pp_{\Sigma(\Gamma^*)}$ given by the equations
$$s_{j+1}t_{j}=s_{j}t_{j-1},\text{ for } j=1,\dots,k.$$
\end{thm}

\begin{pf} Recall from Section~\ref{s:t} that the simplicial toric variety $\ps$ is the geometric
quotient of $\CC^n\setminus{\bf V}(B(\Sigma))$ by $G(\Sigma)$.
Similarly,
$$\pp_{\Sigma(\Gamma^*)}\cong\biggl(\CC^{n+k}\setminus{\bf
V}(B(\Sigma(\Gamma^*))\biggr)/G(\Sigma(\Gamma^*)).$$ We need to
show that the map given by (\ref{e:map}) is well defined, i.e.,
given $(x_1,\dots,x_n)\in G(\Sigma)$, then
$\phi_{\Gamma^*}(x_1,\dots,x_n)\in G(\Sigma(\Gamma^*))$, and,
also, the complement of  ${\bf V}(B(\Sigma))$ maps to the
complement of ${\bf V}(B(\Sigma(\Gamma^*))$. Take
$(x_1,\dots,x_n)$ such that $\prod_{i=1}^n x_i^{\langle
m,e_i\rangle}=1$, for all $m\in M={\rm Hom}(N,\ZZ)$. Then we must
check that $\phi_{\Gamma^*}(x_1,\dots,x_n)$ satisfies the
equations
 \begin{equation}\label{e:g}
 \prod_{e_i\notin\Gamma^*} y_i^{\langle
\tilde m,\tilde e_i\rangle}\prod_{j=1}^{k+1} s_j^{\langle\tilde
m,\tilde e_{s_j}\rangle}t_{j-1}^{\langle\tilde m,\tilde
e_{t_{j-1}}\rangle}=1
 \end{equation}
 for
all $\tilde m\in {\rm Hom}(N\oplus
 \ZZ\{v_1,\dots,v_k\},\ZZ)$. It is
 enough to consider two cases: $\tilde m$ vanishes on
 $v_1,\dots,v_k$ or $N$.
 In the first case, substituting (\ref{e:map}) into the last
 equation, we get
\begin{align*}&
\biggl( \prod_{e_i\notin\Gamma^*} x_i^{\langle \tilde
m,e_i\rangle} \biggr)
 {(x_{i_0}\cdots x_{i_k})}^{\langle\tilde m,e_{i_0}\rangle}\biggl(\prod_{j=1}^{k+1}{(x_{i_{j}}\cdots
x_{i_{k+1}})}^{\langle\tilde
m,\frac{e_{i_{k+1}}-e_{i_0}}{k+1}\rangle}\biggr)x_{i_{k+1}}^{\langle\tilde
m,e_{i_0}\rangle}\hfill
\\&=
\prod_{e_i\notin\Gamma^*} x_i^{\langle \tilde m,e_i\rangle}
\prod_{j=0}^{k+1} x_{i_j}^{\langle \tilde m,
e_{i_0}+\frac{j}{k+1}( e_{i_{k+1}}-e_{i_0})\rangle}=
\prod_{e_i\notin\Gamma^*} x_i^{\langle \tilde
m,e_i\rangle}\prod_{j=0}^{k+1} x_{i_j}^{\langle \tilde
m,e_{i_j}\rangle}
\\&
=\prod_{i=1}^n x_i^{\langle \tilde m,e_i\rangle}=1,
 \end{align*}
 because $\tilde m$ is a homomorphism from $N$ to $\ZZ$. For a
 later use notice that this calculation also shows: if $\phi_{\Gamma^*}(x_1,\dots,x_n)\in
 G(\Sigma(\Gamma^*))$, then $(x_1,\dots,x_n)\in G(\Sigma)$.
In the second case, when $\tilde m$ is zero on $N$, the left-hand
side of (\ref{e:g}) becomes
\begin{align*}
 &\biggl(\prod_{j=1}^{k}
 {(x_{i_0}\cdots x_{i_{j-1}})}^{\langle\tilde m,v_j\rangle}\biggr){(x_{i_0}\cdots x_{i_k}\cdot
 {x_{i_{k+1}}})}^{\langle \tilde m,
 -v_1-\cdots-v_k \rangle}\biggl( \prod_{j=1}^k
{(x_{i_{j}}\cdots x_{i_{k+1}})}^{\langle\tilde
m,v_j\rangle}\biggr)
\\&=
\biggl(\prod_{j=1}^{k}
 {(x_{i_0}\cdots x_{i_{k+1}})}^{\langle\tilde m,v_j\rangle}\biggr){(x_{i_0}\cdots
 {x_{i_{k+1}}})}^{\langle \tilde m,
 -v_1-\cdots-v_k \rangle}
=1.
 \end{align*}

Now, suppose that $(x_1,\dots,x_n)\in\CC^n\setminus{\bf
V}(B(\Sigma))$, then $\prod_{e_i\not\in\sigma} x_i \ne0$ for some
$\sigma\in\Sigma$. If $\sigma$ does not intersect the edge
$\Gamma^*$, then $\sigma\in\Sigma(\Gamma^*)$, and, by the
construction of $\phi_{\Gamma^*}$, all variables corresponding to
the generators in $\Sigma(\Gamma^*)$ outside $\sigma$ don't
vanish, whence so does the product of them. If $\sigma$ intersects
$\Gamma^*$, then, by the construction of $\phi_{\Gamma^*}$ and  of
the corresponding cone $\sigma'\in\Sigma(\Gamma^*)$ of type
(ii)-(iii), the product of the variables for generators outside
$\sigma'$ will not vanish, so that $\phi_{\Gamma^*}(x)$ belongs to
the complement of ${\bf V}(B(\Sigma(\Gamma^*)$ in  $\CC^{n+k}$
 This finishes the proof that the map
$\phi_{\Gamma^*}$ is well defined.

  To prove  injectivity of  $\phi_{\Gamma^*}$  one needs to
recover all $x_i$ from given  $s_1,\dots,s_{k+1}$,
$t_0,\dots,t_{k}$, and $y_i$, accounting for the group actions
which represent the toric varieties as geometric quotients. In
other words, the ``upstairs'' map
$$\phi_{\Gamma^*}:\CC^n\setminus{\bf V}(B(\Sigma)) \rightarrow
\CC^{n+k}\setminus{\bf V}(B(\Sigma(\Gamma^*))$$ must be
one-to-one, and $\phi_{\Gamma^*}(G(\Sigma))=G(\Sigma(\Gamma^*))$.
Note that
\begin{multline}\label{e:eq}
x_{i_j}=\frac{s_{j+1}}{s_j}=\frac{t_{j-1}}{t_j},\text{ for
}j=1,\dots,k,\quad x_{i_0}=s_1,\, x_{i_{k+1}}=t_k,\text{ and }
x_i=y_i
\end{multline}
 for $ e_i\notin\Gamma^*$.
  Hence, $\phi_{\Gamma^*}$ is
injective on $\CC^n\setminus{\bf V}(B(\Sigma))$,  because $s_j$
and $t_j$ don't vanish simultaneously by the properties of toric
varieties: they correspond to the generators $\tilde e_{s_j}$ and
$\tilde e_{t_j}$ which do not lie in the same cone of the fan
$\Sigma(\Gamma^*)$, by the construction. To show the equality take
a point from $G(\Sigma(\Gamma^*))$ with homogeneous coordinates
$s_1,\dots,s_{k+1}$, $t_0,\dots,t_{k}$ and $y_i$, for
$e_i\notin\Gamma^*$. Then, taking $\tilde m$ in  (\ref{e:g}) such
that $\tilde m$ vanishes on $N$ and  $\langle\tilde
m,v_{j}\rangle=-1$, $\langle\tilde m,v_{j+1}\rangle=1$,
$\langle\tilde m,v_l\rangle=0$, for $l\ne j,j+1$, we get the
equations $${s_j^{-1}}{s_{j+1}}t_{j-1}^{-1}{t_j}=1,\text{ for
}j=1,\dots,k-1.$$
 Similarly, if $\tilde m$ is zero on $N$ and on $v_1\dots,v_{k-1}$, and
 $\langle\tilde m,v_k\rangle=-1$,
 then (\ref{e:g}) becomes
$ {s_k^{-1}}{s_{k+1}}t_{k-1}^{-1}{t_k}=1$. Because of
(\ref{e:eq}), these equations show that every point in
$G(\Sigma(\Gamma^*))$ is in the image of $\phi_{\Gamma^*}$.
 Combining this with an earlier note in the proof that the intersection of the
  image of $\phi_{\Gamma^*}$ with $G(\Sigma(\Gamma^*))$ is already
  in the image of $G(\Sigma)$,
  we get
 the equality $\phi_{\Gamma^*}(G(\Sigma))=G(\Sigma(\Gamma^*))$.

Next, we need to verify that $\phi_{\Gamma^*}(\ps)$ coincides with
the complete intersection $s_{j+1}t_{j}=s_{j}t_{j-1},$ for
$j=1,\dots,k.$ One direction is easy by checking that these
equations hold for the points in (\ref{e:map}).
 Conversely, suppose that a
point in $\pp_{\Sigma(\Gamma^*)}$ with homogeneous coordinates
$s_1,\dots,s_{k+1}$, $t_0,\dots,t_{k}$, and $y_i$ satisfies the
above equations. Then we can determine $(x_1,\dots,x_n)$ uniquely
using the equations (\ref{e:eq}). The problem is whether
$\phi_{\Gamma^*}(x_1,\dots,x_n)$ gives us back the original
$s_1,\dots,s_{k+1}$, $t_0,\dots,t_{k}$, and $y_i$. If all $s$ and
$t$ are not equal to zero, then the previous paragraph has already
showed this. So, let us assume that one of $t$, for instance, is
equal to zero and let $k+1\ge l\ge1$ be the minimal integer for
which $t_l\ne0$ and $t_0=\cdots=t_{l-1}=0$. Then notice that if
$t_{j-1}\ne0$ then $t_j\ne0$ as well. Indeed, we have equation
$s_{j+1}t_{j}=s_{j}t_{j-1}$.  Since $s_j$ and $t_j$ don't vanish
simultaneously, we get a contradiction if $t_{j-1}\ne0$ and
$t_j=0$. So, we can assume that $t_0=\cdots=t_{l-1}=0$ and
$t_j\ne0$ for $j\ge l$. Also, then $s_j\ne0$ for $j\le l-1$.
Equation $s_{l+1}t_{l}=s_{l}t_{l-1}$ implies that $s_{l+1}=0$
since $t_l\ne 0$ and $t_{l-1}=0$. Similarly to what we did with
$t_j$, we conclude that $s_{j}=0$ for $j\ge l+1$. Now, plugging
$x_{i_j}={t_{j-1}}/{t_j}$, for $j=l,\dots,k$, and
$x_{i_{k+1}}=t_k$ into $t_j=x_{i_{j+1}}x_{i_{j+2}}\dots
x_{i_k}x_{i_{k+1}}$ we recover $t_j$ for $j\ge l-1$. Doing the
same for $s_j=x_{i_0}x_{i_1}\dots x_{i_{j-2}}x_{i_{j-1}}$ with
$x_{i_j}={s_{j+1}}/{s_j}$ and $x_{i_0}=s_1$ we get $s_j$ for $j\le
l$. Since $x_{i_l}=t_{l-1}/t_l=0/t_l=0$, equations (\ref{e:map})
give us $s_j=0$ for $j\ge l+1$ and $t_j=0$ for $j\le l-1$. Thus,
each point in the complete intersection is in the image
 $\phi_{\Gamma^*}(\ps)$.

To finish the proof  we are left to show that the complete
intersection $s_{j+1}t_{j}=s_{j}t_{j-1},$ for $j=1,\dots,k,$ is
quasismooth.  Quasismooth complete intersections in toric
varieties were introduced in \cite[Definition~1.1]{m0} and  are
equivalent to $V$-submanifolds (suborbifolds; see \cite{bc}) of
the ambient toric variety. Recall that a complete intersection
defined by the homogeneous polynomials $f_1,\dots,f_s$ is {\it
quasismooth} if the  Jacobian matrix $\bigl(\frac{\partial
f_j}{\partial x_i}(x)\bigr)_{i,j}$ has maximal rank on the
intersection. In our case the Jacobian matrix is $$ \left(
 \begin{array}{cccccccccccccc}
 t_0& -t_1 &0&0&\cdots&0& s_1& -s_2 &0&0&\cdots&0&0&\cdots
\\0& t_1&-t_2&0&\cdots&0&0& s_2&-s_3&0&\cdots&0&0&\cdots
\\0&0& t_2&-t_3&\cdots&0&0&0& s_3&-s_4&\cdots&0&0&\cdots
\\\vdots&&&\ddots&\ddots&\vdots&\vdots&&&\ddots&\ddots&\vdots&0&\cdots
\\0&0&0&\cdots&t_{k-1}&-t_k&0&0&0&\cdots& s_{k}&-s_{k+1}&0&\cdots\\
\end{array}\right)$$
where rows correspond to the  partial derivatives of the defining
polynomials $s_{j+1}t_{j}-s_{j}t_{j-1}$,  $j=1,\dots,k$, with
respect to the variables $s_1,\dots,s_{k+1},t_0,\dots,t_{k}$ and
$y_i$ for $e_i\notin \Gamma^*$. As in the previous paragraph, we
can assume that a point on the complete intersection will have
coordinates  $s_j\ne0$ for $j\le l-1$ and $t_j\ne0$ for $j\ge l$
for some $0\le l\le k+1$. Then it is easy to see that the above
matrix will have the maximal rank equal to $k$.
 The proof of the theorem is complete.
\end{pf}

\begin{rem}\label{r:comp}
If $X$ is an anticanonical  hypersurface
defined by a homogeneous polynomial $$f(x)=\sum_{m\in\Delta\cap M}a_m
\prod_{i=1}^n x_i^{1+\langle m,e_i\rangle}$$
in the toric variety $\ps$ satisfying the condition of
Theorem~\ref{t:emb}, then by the map $\phi_{\Gamma^*}$ the variety $X$ is isomorphic to
a complete intersection in $\pp_{\Sigma(\Gamma^*)}$ given by
$$s_{j+1}t_{j}=s_{j}t_{j-1},\text{ for } j=1,\dots,k,$$
and
\begin{equation}\label{e:pol}
\sum_{m\in\Delta\cap M}a_m X_m
\prod_{e_i\notin\Gamma^*} y_i^{1+\langle m,e_i\rangle}=0,
\end{equation}
 where $$X_m=  \left\{
 \begin{matrix}
(s_1t_0)^{1+\langle m,e_{i_0}\rangle}
(t_0\cdots t_k)^{\langle m,e_{i_1}-e_{i_0}\rangle}  & \mbox{ if } &\langle m,e_{i_1}-e_{i_0}\rangle>0
\\
(s_1t_0)^{1+\langle m,e_{i_{k+1}}\rangle}
(s_1\cdots s_{k+1})^{\langle m,e_{i_0}-e_{i_1}\rangle}
   & \mbox{ if } &\langle m,e_{i_1}-e_{i_0}\rangle<0.
\end{matrix}\right.$$
\end{rem}

When a $V$-manifold (orbifold) is embedded into another
$V$-manifold it need not to be a $V$-submanifold (see Remark 3.3
in \cite{bc}). Magically, again, the complete intersection which
represents the $V$-manifold $X$ in the toric variety
$\pp_{\Sigma(\Gamma^*)}$ is quasismooth, i.e., $V$-submanifold, if
$X$ were generic enough (at least nondegenerate).

\begin{pr} If $X$ is a semiample anticanonical nondegenerate hypersurface defined by a
homogeneous polynomial $$f(x)=\sum_{m\in\Delta\cap M}a_m
\prod_{i=1}^n x_i^{1+\langle m,e_i\rangle}$$ in a complete
simplicial toric variety $\ps$ as in Theorem~\ref{t:emb}, then the
complete intersection variety $X$ described in Remark~\ref{r:comp}
 is quasismooth in the toric variety
 $\pp_{\Sigma(\Gamma^*)}$.
 \end{pr}

\begin{pf} To show that the complete intersection is
quasismooth we need to check if the Jacobian matrix consisting of
partial derivatives of defining polynomials has maximal possible
rank on the  complete intersection. In our case it will
essentially suffice to look at the part of this matrix
corresponding to the derivatives with respect to
$s_1,\dots,s_{k+1},t_0,\dots,t_{k}$: $$\left(
\begin{array}{cccccccccccccc}
 t_0& -t_1 &0&0&\cdots&0& s_1& -s_2 &0&0&\cdots&0
\\0& t_1&-t_2&0&\cdots&0&0& s_2&-s_3&0&\cdots&0
\\0&0& t_2&-t_3&\cdots&0&0&0& s_3&-s_4&\cdots&0
\\\vdots&&&\ddots&\ddots&\vdots&\vdots&&&\ddots&\ddots&\vdots
\\0&0&0&\cdots&t_{k-1}&-t_k&0&0&0&\cdots& s_{k}&-s_{k+1}
\\
\frac{\partial {\tilde f}}{\partial s_1}&& &\cdots&&
\frac{\partial \tilde f}{\partial s_{k+1}}&\frac{\partial \tilde
f}{\partial t_0}&& &\cdots&& \frac{\partial \tilde f}{\partial
t_{k}}\\
\end{array}\right)$$
where $\tilde f$  denotes the polynomial on the right hand side of
the equation (\ref{e:pol}). What is missing in the above matrix
compared to the whole Jacobian matrix are the partial derivatives
with respect to $y_i$ for  $e_i\notin \Gamma^*$.

 First, note that
$\partial \tilde f/\partial y_i$ coincides with $\partial
f/\partial x_i$ for $e_i\notin \Gamma^*$ after transformation
(\ref{e:map}). So, if one of such $\partial f/\partial x_i$ does
not vanish then the Jacobian matrix will have the maximal rank
equal to $k+1$, since the above matrix without the last row is
already shown to have the maximal rank on the intersection in the
proof of Theorem~\ref{t:emb}.
 Now, since $X$ is
nondegenerate, polynomials   $x_1\partial f/\partial
x_1,\dots,x_n\partial f/\partial x_n$ do not vanish simultaneously
on $\ps$. Hence, we can assume that one of the following
polynomials $x_{i_0}\partial f/\partial
x_{i_0},\dots,x_{i_{k+1}}\partial f/\partial x_{i_{k+1}}$ does not
vanish on $\ps$. Then, since every $x_{i_l}\partial f/\partial
x_{i_l}$, for $l=1,\dots,k$, can be expressed via $x_{i_0}\partial
f/\partial x_{i_0}$ and $x_{i_{k+1}}\partial f/\partial
x_{i_{k+1}}$ by Euler formulas (see Definition~3.9 in \cite{bc}),
the polynomial $x_{i_0}\partial f/\partial x_{i_0}$ or
$x_{i_{k+1}}\partial f/\partial x_{i_{k+1}}$ is not zero on  $X$.

Second, note that $x_{i_0}\partial f/\partial x_{i_0}$ coincides
with the polynomial $$s_1\frac{\partial {\tilde f}}{\partial
s_1}+\cdots+s_{k+1}\frac{\partial \tilde f}{\partial s_{k+1}}$$
after the substitution (\ref{e:map}). (This is rather a
straightforward calculation that we leave to the reader.)
Similarly, $x_{i_{k+1}}\partial f/\partial x_{i_{k+1}}$
corresponds to $$t_0\frac{\partial {\tilde f}}{\partial
t_0}+\cdots+t_{k}\frac{\partial \tilde f}{\partial t_{k}}.$$ Now,
as in the proof of the above theorem,  a point on the complete
intersection in $\pp_{\Sigma(\Gamma^*)}$ will either have all
coordinates $s$ and $t$ not equal to zero, or  we can assume that
there is $l$ between $0$ and $k+1$ such that $t_j\ne0$, for $j\ge
l$, $t_0=\cdots=t_{l-1}=0$, $s_j\ne0$, for $j\le l-1$, and
$s_{l+1}=\cdots s_{k+1}=0$. In the latter case, it will suffice to
have that one of the partial derivatives $$\frac{\partial {\tilde
f}}{\partial s_1},\dots,\frac{\partial \tilde f}{\partial
s_{l}},\frac{\partial \tilde f}{\partial t_l},\dots \frac{\partial
\tilde f}{\partial t_{k}}$$ in the last row of the matrix does not
vanish. But if they do vanish, then both  $x_{i_0}\partial
f/\partial x_{i_0}$ and $x_{i_{k+1}}\partial f/\partial
x_{i_{k+1}}$ will be zero by the above representation in terms of
the sums, which would give a contradiction with nondegeneracy of
$X$. If all $s$ and $t$ do not vanish, then multiply each of
column of the above matrix, corresponding to the partial
derivative with respect to a variable $s$ or $t$, by that
variable. This procedure will not change the rank of the matrix.
The left half of the new matrix will be
 $$\left(
\begin{array}{cccccccccccccc}
 s_1t_0& -s_2t_1 &0&0&\cdots&0
\\0& s_2t_1&-s_3t_2&0&\cdots&0
\\0&0& s_3t_2&-s_4t_3&\cdots&0
\\\vdots&&&\ddots&\ddots&\vdots
\\0&0&0&\cdots&s_kt_{k-1}&-s_{k+1}t_k
\\
s_1\frac{\partial {\tilde f}}{\partial s_1}&& &\cdots&&
s_{k+1}\frac{\partial \tilde f}{\partial s_{k+1}} \\
\end{array}\right)$$

Since $s_{j+1}t_{j}=s_{j}t_{j-1}$,  for $j=1,\dots,k$, adding in
succession to the second column the first one, then to the third
the second one, and so on, will produce a new matrix: $$\left(
\begin{array}{cccccccccccccc}
 s_1t_0&  0&0&\cdots&&0
\\0& s_1t_0& 0  &\cdots&&0
\\0&0& s_1t_0 &\cdots&&0
\\\vdots&&&\ddots&&\vdots
\\0& &\cdots& &s_1t_0&0
\\
 &&\cdots&&&
* \\
\end{array}\right),$$
where $s_1t_0\ne0$ on the diagonal and $*$ is the polynomial
$s_1\frac{\partial {\tilde f}}{\partial
s_1}+\cdots+s_{k+1}\frac{\partial \tilde f}{\partial s_{k+1}}$,
which by symmetry we can assume doesn't vanish. Hence, the
original Jacobian matrix has the maximal rank.
\end{pf}

The second main result shows that deformations of the toric
variety $\ps$ and its Calabi-Yau hypersurface are also realized as
complete intersections inside another toric variety. The proof of the
theorem also shows a relation of these deformations to the
deformations constructed in \cite{m3}.

\begin{thm}\label{t:2ndmain}
Let $X$ be a semiample anticanonical  hypersurface, defined by a
homogeneous polynomial $$f(x)=\sum_{m\in\Delta\cap M}a_m
\prod_{i=1}^n x_i^{1+\langle m,e_i\rangle}$$ in $S(\Sigma)$, in a
complete simplicial toric variety $\ps$ as in Theorem~\ref{t:emb}.
Then the following family of complete intersections in
$\pp_{\Sigma(\Gamma^*)}$ given by $$
s_{j+1}t_{j}=s_{j}t_{j-1}+\sum_{u\in{\rm
int}(\Gamma)}\lambda_{u,j}\prod_{e_i\notin\Gamma^*} y_i^{\langle
u,e_i\rangle},\quad\text{ for } j=1,\dots,k,$$ and equation
(\ref{e:pol})
 are  non-polynomial deformations  of $X$, corresponding
 to
 $$\bigoplus\begin{Sb}
e_l\in{{\rm int}(\Gamma^*)}\end{Sb}(S/\langle
x_l\rangle)_{\beta_1^\sigma}\subset H^1(X,{\cal T}_X)$$ from Proposition~\ref{t:chiral} by the
Kodaira-Spencer map. Polynomial deformations of $X$  are parametrized by the coefficients $a_m$ in (\ref{e:pol}).
\end{thm}

\begin{pf} Using a different method in \cite{m3}, we already  constructed 1-parameter non-polynomial
deformations of the Calabi-Yau hypersurface $X$, which are induced
by the deformations of the toric variety $\ps$. By the
Kodaira-Spencer map, the corresponding infinitesimal deformations
in $H^1(X,{\cal T}_X)$  have also been found in Section 3 of that
paper. The essential work left is to show that those deformations of the toric variety coincide
with the following 1-parameter subfamilies of deformations of
$\ps$, which are complete intersections in $\pp_{\Sigma(\Gamma^*)}$:
\begin{equation}\label{e:1stfam}
s_{j+1}t_{j}=s_{j}t_{j-1},\text{ for } j\ne l\,\text{ and }\,
s_{l+1}t_{l}=s_{l}t_{l-1}+\lambda\prod_{e_i\notin\Gamma^*}
y_i^{\langle u,e_i\rangle},
\end{equation}
where $u\in{\rm int}(\Gamma)$.
After this, we only need to check that everything carries over for the hypersurface.

 Recall the construction from
\cite{m3}. Let $U^l_0=\{x\in\ps: \prod_{k>j}x_{l_k}\ne0\} $ and
$U^l_1=\{x\in\ps: \prod_{k<j}x_{l_k}\ne0\}$ be  the open cover of
the toric variety $\ps$. The deformation of $\ps$ is defined by
regluing $U^l_0$ and $U^l_1$ along their intersection $U^l_0\cap
U^l_1$. More precisely, the 1-parameter family ${\cal P}$ is
obtained by gluing
 $U^l_0\times{\AA}^1$ and
$U^l_1\times{\AA}^1$  along $U^l_0\cap U^l_1\times{\AA}^1$ by the
identification:
\begin{equation}\label{e:id}
(x_1,\dots,x_l,\dots,x_n,\lambda)\mapsto (x_1,\dots,x_l-\lambda
x_l\prod_{i=1}^n x_i^{\langle u,e_i\rangle},\dots,x_n,\lambda).
\end{equation}
(Note: instead of $2\lambda$ in \cite{m3} we have $\lambda$ here
for simplification.)

The complete intersections (\ref{e:1stfam}) can also be put into a
family over ${\AA}^1$:
 $${\cal P}'=\{(y,\lambda)\in\pp_{\Sigma(\Gamma^*)}\times\AA^1:\,(\ref{e:1stfam}) \text{ holds for }
 (y,\lambda)\}.$$
 Let us show that the families ${\cal P}$ and ${\cal P}'$ are
 isomorphic.
The second family ${\cal P}'$ admits a covering by two open
subsets: $\tilde U^l_0$ and $\tilde U^l_1$ defined by  $t_l\ne0$
and $s_l\ne0$, respectively. Consider the following isomorphism
$U^l_0\times\AA^1\cong\tilde U^l_0$, defined by composing
$\phi_{\Gamma^*}\times{\rm id}$ with the map on $\tilde U^l_0$
which changes only the variables $$s_j\mapsto
s_j-\frac{\lambda}{t_{j-1}}\prod_{e_i\notin\Gamma^*} y_i^{\langle
u,e_i\rangle}$$ for $j\ge l+1$. Note that these $t_{j-1}\ne0$
since $t_l\ne0$, as in the proof of Theorem~\ref{t:emb}.
Similarly, $U^l_1\times\AA^1\cong\tilde U^l_1$,which is a
composition of $\phi_{\Gamma^*}\times{\rm id}$ with the map on
$\tilde U^l_1$ which changes $$t_j\mapsto
t_j+\frac{\lambda}{s_{j+1}}\prod_{e_i\notin\Gamma^*} y_i^{\langle
u,e_i\rangle}$$ for $j\le l-1$. Since  $u\in{\rm int}(\Gamma)$, by
the properties of the dual faces $\langle u,e_i\rangle=-1$ for
$e_i\in\Gamma^*$. It is straightforward to check that the isomorphim
$U^l_0\times\AA^1\cong\tilde U^l_0$ is the composition of
(\ref{e:id}) with the above defined $U^l_1\times\AA^1\cong\tilde
U^l_1$.

The 1-parameter family ${\cal X}\subset{\cal P}$ corresponding to
the deformation of the Calabi-Yau hypersurface $X$, determined by
 $$f(x)=\sum_{m\in\Delta\cap M}a_m
\prod_{i=1}^n x_i^{1+\langle m,e_i\rangle},$$ is given by the
hypersurfaces in $U^l_\varepsilon\times\AA^1$, $\varepsilon=0,1$,
corresponding to $$f^\varepsilon_\lambda(x)=\sum_{m\in\Delta}a_m
\prod_{i=1}^n x_i^{-1+\langle m,e_i\rangle}\biggl(1+
((-1)^\varepsilon+c_m)\frac{\lambda}{2}\prod_{i=1}^n
 x_i^{\langle u,e_i\rangle}\biggr)^{-1+\langle m,e_l\rangle},$$
where
 $$c_m= \left\{
 \begin{matrix}
 1 & \mbox{ if } &\langle m,e_{l_0}-e_l\rangle)>0
\\
 0 & \mbox{ if } &\langle m,e_{l_0}-e_l\rangle)=0
 \\
  -1 & \mbox{ if } &\langle m,e_{l_0}-e_l\rangle)<0.
\end{matrix}\right.$$
One can check that these hypersurfaces are compatible with respect
to the identification (\ref{e:id}). Moreover, polynomials  $f^\varepsilon_\lambda(x)$ will both
transfer to the polynomial
$$\sum_{m\in\Delta\cap M}a_m X_m
\prod_{e_i\notin\Gamma^*} y_i^{1+\langle m,e_i\rangle}$$
 in
(\ref{e:pol}) by the already defined isomorphisms
$U^l_\varepsilon\times\AA^1\cong\tilde U^l_\varepsilon$. Hence,
this  deformation of $X$ is isomorphic to the intersection of
(\ref{e:1stfam}) with (\ref{e:pol}).

The fact that the polynomial deformations of $X$ correspond to
deformations  of  $f$ is shown in \cite[Section~1]{m3}.
\end{pf}

\begin{rem}
In the above theorem, the non-polynomial deformations for a
minimal Calabi-Yau hypersurface $X$ correspond to the term
$l^*(\Gamma)l^*(\Gamma^*)$ in Proposition~\ref{p:dim}.
\end{rem}

\begin{rem}
In the Batyrev-Borisov mirror symmetry construction (see
\cite{bb}), the Calabi-Yau complete intersections are nef. In our
construction, we can assume that the ambient toric variety
$\pp_{\Sigma(\Gamma^*)}$  is complete (by compactifying it if
necessary). However, in most cases, the complete intersection in
Theorem~\ref{t:2ndmain} will not be nef. Here is the simplest
example of such a Calabi-Yau complete intersection. Let $\Sigma$
be the fan with cones formed by a triangulation of the reflexive
polytope whose vertices are
$$\{(-2,-2,-3,-2),(1,0,0,0),(0,1,0,0),(0,0,1,0),(0,0,1,2)\}.$$
Except for the origin there are three more lattice points
$(0,0,1,1)$, $(-1,-1,-1,-1)$ and $(-1,-1,-1,0)$ in this polytope
that lie in the middle of the edges of the two-dimensional face
with vertices at $(0,0,1,0)$, $(0,0,1,2)$ and $(-2,-2,-3,-2)$. We
can apply our construction to the edge  between
$e_{i_0}=(0,0,1,0)$ and $e_{i_2}=(0,0,1,2)$ with
$e_{i_1}=(0,0,1,1)$ as the only interior lattice point. By the
formulas $(\ref{e:generators})$, we get a new fan
$\Sigma(\Gamma^*)$ with the set of minimal integral generators at
the points
\begin{multline*}
(-2,-2,-3,-2,0),(1,0,0,0,0),(0,1,0,0,0),(-1,-1,-1,-1,0),\\
(-1,-1,-1,0,0),\,v=(0,0,0,0,1),\, e_{i_0}-v=(0,0,1,0,-1),\\
v+\frac{e_{i_2}-e_{i_0}}{2}=
(0,0,0,1,1),\,e_{i_0}-v+\frac{e_{i_2}-e_{i_0}}{2}=(0,0,1,1,-1).
\end{multline*}

However, the convex hull of these points in $\ZZ^5$ is not a
reflexive polytope. On the other hand, if the hypersurfaces
defining the complete intersection in $\pp_{\Sigma(\Gamma^*)}$
were nef, then the sum of the degrees of these hypersurfaces was
also nef and anticanonical, by the construction. But in this case,
by the uniqueness result of \cite[Proposition~1.2]{m1} and
\cite[Theorem~4.1.9]{b}, the variety $\pp_{\Sigma(\Gamma^*)}$ must
be a partial crepant resolution of a Fano toric variety, in which
case, the convex hull of the generators of $\Sigma(\Gamma^*)$ form
a reflexive polytope. Hence, we get a contradiction with the
assumption that the complete intersection was nef. This example is
very important. It indicates that we have to study non-nef
Calabi-Yau complete intersections in toric varieties.
\end{rem}

{\it Acknowledgement.} The  author   would like to thank the
Max-Planck Institute for Mathematics in Bonn for an opportunity to
finish this paper.

\end{document}